\documentclass[12pt]{amsart}
\usepackage{amsmath,amssymb,color,placeins,cite,array,multirow,colonequals,verbatim}
\usepackage{etoolbox}
\apptocmd{\sloppy}{\hbadness 10000\relax}{}{}

\numberwithin{equation}{section}

\newtheorem{thm}[equation]{Theorem}

\newtheorem{prop}[equation]{Proposition}
\newtheorem{lemma}[equation]{Lemma}

\theoremstyle{definition}
\newtheorem{rmk}[equation]{Remark}
\newtheorem{notation}[equation]{Notation}
\newtheorem{defn}[equation]{Definition}

\newcommand{\F}{\mathbb{F}}
\newcommand{\bP}{\mathbb{P}}

\DeclareMathOperator{\charp}{char}

\DeclareMathOperator{\Tr}{Tr}

\usepackage[colorlinks,pagebackref,pdftex,bookmarks=false]{hyperref}

\begin{document}

\title{Some classes of permutation pentanomials}

\author{Zhiguo Ding}
\address{
  School of Mathematics and Statistics, 
  Central South University, 
  Changsha, 410075 China
}
\email{ding8191@csu.edu.cn}

\author{Michael E. Zieve}
\address{
  Department of Mathematics,
  University of Michigan,
  530 Church Street,
  Ann Arbor, MI 48109-1043 USA
}
\email{zieve@umich.edu}
\urladdr{https://dept.math.lsa.umich.edu/$\sim$zieve/}

\date{\today}

\begin{abstract}
For each prime $p\ne 3$ and each power $q=p^k$, we present two large classes of permutation polynomials over $\F_{q^2}$ of the form $X^r B(X^{q-1})$ which have at most five terms, where $B(X)$ is a polynomial with coefficients in $\{1,-1\}$. 
The special case $p=2$ of our results comprises a vast generalization of $76$ recent results and conjectures in the literature. In case $p>2$, no instances of our permutation polynomials have appeared in the literature, and the construction of such polynomials had been posed as an open problem. Our proofs are short and involve no computations, in contrast to the proofs of many of the
special cases of our results which were published previously.
\end{abstract}

\thanks{
The first author was supported in part by the Natural Science Foundation of Hunan Province of China (No.\ 2020JJ4164).
The second author was supported in part by Simons Travel Grant MPS-TSM-00007931.
}

\maketitle


\section{Introduction}

A polynomial $f(X)\in\F_q[X]$ is called a \emph{permutation polynomial} if the function $\alpha\mapsto f(\alpha)$ permutes $\F_q$. In recent years, many authors have constructed classes of permutation polynomials having at most five terms and having coefficients in $\{1,-1\}$. In particular, many of these are permutation polynomials over $\F_{q^2}$ with the form $X^r B(X^{q-1})$. 
In this paper we construct two large classes of such permutation polynomials, which subsume $76$ recent results and conjectures in the literature (from $20$ papers by $47$ authors). Notably, all of these previous results and conjectures restrict to the case that $\F_q$ has characteristic $2$. 
In characteristic $2$, our results comprise vast generalizations of all these previous results and conjectures; while at the same time, our results show that characteristic $2$ is not special in this context, since our results simultaneously yield permutation polynomials of the same form in every characteristic other than $3$. 
Our permutation polynomials in characteristic greater than $3$ seem to be completely new, in the sense that we could not find any instances of these permutation polynomials in the literature. Moreover, our permutation polynomials in odd characteristic resolve Open Problem 2 from \cite[\S 5]{KX} (we will resolve Open Problem 1 elsewhere).
%
%
Our main results are as follows.

\begin{thm}\label{1}
Let $p$ be prime with $p\ne 3$, let $q=p^k$ where $k$ is a positive integer, and pick $Q,R,S\in\{p^i:i\ge 0\}$. Pick an order-$3$ element $\omega\in\F_{q^2}^*$, and let $r$ be a positive integer such that $r\equiv Q+R+S\pmod{q+1}$. Write
\begin{align*}
C_1(X)\colonequals &\,
(\omega^{Q+R+S}-\omega) X^{Q+R+S}
+ (\omega^{Q+R}-\omega^{S+1}) X^{Q+R} \\
&\,+ (\omega^{Q+S}-\omega^{R+1}) X^{Q+S}
+ (\omega^{R+S}-\omega^{Q+1}) X^{R+S} \\
&\,+ (\omega^Q-\omega^{R+S+1}) X^Q
+ (\omega^R-\omega^{Q+S+1}) X^R \\
&\,+ (\omega^S-\omega^{Q+R+1}) X^S
+ (1-\omega^{Q+R+S+1})
\end{align*}
and $C_2(X)\colonequals X^{Q+R+S} C_1(1/X)$.
Let $\beta=\omega^{Q+R}-\omega^{S+1}$ if $Q+R+S\equiv 1\pmod 3$, and $\beta=\omega^{Q+R+S}-\omega$ otherwise, and write $B_z(X)\colonequals C_z(X)/\beta$ for each $z\in\{1,2\}$. 
Then $B_z(X)$ has at most five terms, its coefficients are in\/ $\F_p$, and if the integers $Q,R,S$ are pairwise distinct then the coefficients of $B_z(X)$ are in $\{1,-1\}$. 
Let $e$ be the unique integer in $\{1,-1\}$ with $q\equiv e\pmod 3$. Then $f(X)\colonequals X^r B_z(X^{q-1})$ permutes\/ $\F_{q^2}$ if and only if $\gcd(r,q-1)=1=\gcd(Q+R+S,q+e)$.
\end{thm}

In Table~\ref{tab1} we give simpler expressions for the polynomials $B_z(X)$ in Theorem~\ref{1}.

\begin{table}[!htbp]
\setlength\extrarowheight{3pt}
\caption{The polynomials $B_z(X)$, where $\sigma=(a,b,c)$ with $a,b,c\in\{1,-1\}$ and $Q\equiv a\pmod 3$, $R\equiv b\pmod 3$, $S\equiv c\pmod 3$}
\label{tab1}
\begin{tabular}{c|l}
\hline
$\sigma$&$B_z(X)$ \\ \hline
\multirow{2}{*}{$(1,1,1)$}&
$\displaystyle{B_1(X)=X^{Q+R+S} - X^Q - X^R - X^S + 1}$ \\
&$\displaystyle{B_2(X)=X^{Q+R+S} - X^{Q+R} - X^{Q+S} - X^{R+S} + 1}$ \\ \hline
\multirow{2}{*}{$(1,1,-1)$}&
$\displaystyle{B_1(X) = X^{Q+R} - X^{Q+S} - X^{R+S} + X^S - 1 }$ \\
&$\displaystyle{B_2(X) = -X^{Q+R+S} + X^{Q+R} - X^Q - X^R + X^S}$ \\ \hline
\multirow{2}{*}{$(1,-1,-1)$}&
$\displaystyle{B_1(X) = X^{Q+R+S} - X^{R+S} - X^Q + X^R + X^S}$ \\
&$\displaystyle{B_2(X) = X^{Q+R} + X^{Q+S} - X^{R+S} - X^Q + 1}$ \\ \hline
\multirow{2}{*}{$(-1,-1,-1)$}&
$\displaystyle{B_1(X) = X^{Q+R+S} - X^{Q+R} - X^{Q+S} - X^{R+S} + 1}$ \\
&$\displaystyle{B_2(X) = X^{Q+R+S} - X^Q - X^R - X^S + 1}$ \\ \hline
\end{tabular}
\end{table}

\begin{rmk}
In Table~\ref{tab1} it is not necessary to consider other triples $\sigma$, since interchanging either $Q$ and $R$ or $Q$ and $S$ has the effect of multiplying $B_z(X)$ by $\pm 1$.
\end{rmk}

%
%

\begin{thm}\label{2}
Let $p$ be prime with $p\ne 3$, let $q=p^k$ where $k$ is a positive integer, and pick $Q,R,S\in\{p^i:i\ge 0\}$. Pick an order-$3$ element $\omega\in\F_{q^2}^*$, and let $r$ be a positive integer such that $r\equiv Q+R+S\pmod{q+1}$. Write
\begin{align*}
C_1(X)\colonequals &\,
(\omega^{Q+S}-\omega^{R+1}) X^{Q+R+S}
+ (\omega^Q-\omega^{R+S+1}) X^{Q+R} \\
&\,+ (\omega^{Q+R+S}-\omega) X^{Q+S}
+ (\omega^S-\omega^{Q+R+1}) X^{R+S} \\
&\,+ (\omega^{Q+R}-\omega^{S+1}) X^Q
+ (1-\omega^{Q+R+S+1}) X^R \\
&\,+ (\omega^{R+S}-\omega^{Q+1}) X^S
+ (\omega^R-\omega^{Q+S+1})
\end{align*}
and $C_2(X)\colonequals X^{Q+R+S} C_1(1/X)$.
Let $\beta=\omega^Q-\omega^{R+S+1}$ if $Q+S\equiv R+1\pmod 3$, and $\beta=\omega^{Q+S}-\omega^{R+1}$ otherwise, and write $B_z(X)\colonequals C_z(X)/\beta$ for each $z\in\{1,2\}$. 
Then $B_z(X)$ has at most five terms, its coefficients are in\/ $\F_p$, and if the integers $Q,R,S$ are pairwise distinct then the coefficients of $B_z(X)$ are in $\{1,-1\}$. Moreover, $f(X)\colonequals X^r B_z(X^{q-1})$ permutes\/ $\F_{q^2}$ if and only if $q\equiv 1\pmod 3$ and $\gcd(r,q-1)=1=\gcd(Q-R+S,q+1)$.
\end{thm}

In Table~\ref{tab2} we give simpler expressions for the polynomials $B_z(X)$ in Theorem~\ref{2}.

\begin{table}[!htbp]
\setlength\extrarowheight{3pt}
\caption{The polynomials $B_z(X)$, where $\sigma=(a,b,c)$ with $a,b,c\in\{1,-1\}$ and $Q\equiv a\pmod 3$, $R\equiv b\pmod 3$, $S\equiv c\pmod 3$}
\label{tab2}
\begin{tabular}{c|l}
\hline
$\sigma$&$B_z(X)$ \\ \hline
\multirow{2}{*}{$(1,1,1)$}&
$\displaystyle{B_1(X)=X^{Q+R} - X^{Q+S} + X^{R+S} - X^R + 1}$ \\
&$\displaystyle{B_2(X)=X^{Q+R+S} - X^{Q+S} + X^Q - X^R + X^S}$ \\ \hline
\multirow{2}{*}{$(1,1,-1)$}&
$\displaystyle{B_1(X) = X^{Q+R+S} - X^{R+S} - X^Q + X^R + X^S}$ \\
&$\displaystyle{B_2(X) = X^{Q+R} + X^{Q+S} - X^{R+S} - X^Q + 1}$ \\ \hline
\multirow{2}{*}{$(1,-1,1)$}&
$\displaystyle{B_1(X) = X^{Q+R+S} - X^Q - X^R - X^S + 1}$ \\
&$\displaystyle{B_2(X) = X^{Q+R+S} - X^{Q+R} - X^{Q+S} - X^{R+S} + 1}$ \\ \hline
\multirow{2}{*}{$(1,-1,-1)$}&
$\displaystyle{B_1(X) = X^{Q+R} - X^{Q+S} - X^{R+S} + X^S - 1}$ \\
&$\displaystyle{B_2(X) = -X^{Q+R+S} + X^{Q+R} - X^Q - X^R + X^S}$ \\ \hline
\multirow{2}{*}{$(-1,1,-1)$} &
$\displaystyle{B_1(X) = X^{Q+R+S} - X^{Q+R} - X^{Q+S} - X^{R+S} + 1}$ \\
&$\displaystyle{B_2(X) = X^{Q+R+S} - X^Q - X^R - X^S + 1}$ \\ \hline
\multirow{2}{*}{$(-1,-1,-1)$}&
$\displaystyle{B_1(X) = X^{Q+R+S} - X^{Q+S} + X^Q - X^R + X^S}$ \\
&$\displaystyle{B_2(X) = X^{Q+R} - X^{Q+S} + X^{R+S} - X^R + 1}$ \\ \hline
\end{tabular}
\end{table}

\begin{rmk}
In Table~\ref{tab2} it is not necessary to consider the triples $\sigma=(-1,\pm 1,1)$, since interchanging $Q$ and $S$ has the effect of multiplying $B_z(X)$ by $\pm 1$.
\end{rmk}

%
%

As we noted previously, many authors have published instances of the permutation polynomials in Theorems~\ref{1} and \ref{2} in case $p=2$. This was done by considering each possibility for $Q=2^a$, $R=2^b$, $S=2^c$ where $a,b,c$ are prescribed nonnegative integers, sometimes along with each choice of a positive integer $r$ such that $r\equiv Q+R+S\pmod{q+1}$. 
Of course, there are infinitely many choices for the quadruple $(a,b,c,r)$, each of which yields permutation polynomials over $\F_{q^2}$ for infinitely many $q=2^k$. Thus, our Theorems~\ref{1} and \ref{2} should be viewed as two four-parameter sequences of results, where each result consists of a collection of permutation polynomials over $\F_{q^2}$ for infinitely many $q$, and where many of the previous theorems in the literature are individual instances of our four-parameter sequence of results. 
We also note that most of the previous results were proved via arguments that are much longer and more computational than ours. Finally, we emphasize that our results apply to every prime $p$ other than $3$, whereas all previous results applied only to $p=2$.  

We now describe the previous results in more detail. In Tables \ref{tabpre}--\ref{tabpre4} we list the previous results and conjectures which are subsumed by our Theorems~\ref{1} and \ref{2}. All of these previous results and conjectures address only the case $p=2$. Each row of each table lists a previous result or conjecture, along with which case(s) of which of our theorems it is subsumed by. 
When there is no entry in a column ``$Q$'', ``$R$'', or ``$S$'', it means that the relevant quantity can be an arbitrary element of $\{2^i:i\ge 0\}$ satisfying the conditions listed in the ``Conditions'' column when such a column exists. When there is no entry in column ``$r$'', it means that $r$ can be any positive integer congruent to $Q+R+S$ (mod $q+1$). In each row of each table, we consider two polynomials to be the same if they are congruent to one another mod $X^{q^2}-X$. 
The column ``e'' in Tables~\ref{tabpre} and \ref{tabpre3} indicates the type of equivalence being used to relate the stated result from the literature to an instance of the stated case of our results: if there is no entry then the equivalence is equality; an entry of ``s'' indicates ``scalar equivalence'', meaning that the two polynomials are obtained from one another by composing on both sides with scalar multiples $\gamma X$ with $\gamma\in\F_{q^2}^*$; an entry of ``m'' indicates ``multiplicative equivalence'', meaning that the two polynomials are obtained from one another by composing on the left with a scalar multiple and on the right with a permutation monomial $\gamma X^s$ with $\gamma\in\F_{q^2}^*$ and $\gcd(s,q^2-1)=1$. 
Our tables contain implicit corrections to three of the cited results: one must add the condition that $k$ is odd to \cite[Table~2(14)]{ZZZ}, one must add the condition $\gcd(3,m)=1$ to \cite[Thm.~4.2]{XCP}, and in \cite[Thm.~2]{KX} one should have $r_A=0$ if and only if $i\equiv j\equiv 1\pmod 2$.

\begin{table}[!htbp]
\setlength\extrarowheight{3pt}
\caption{Previous results and conjectures subsumed by our results, I}
\label{tabpre}
\begin{tabular}{l|c|c|c|c|c|c|c|l}
\hline
Reference&e&Thm.&$B_z$&$Q$&$R$&$S$&$r$&Conditions \\
\hline
\hline
\cite[Thm.~4.2]{ZHF}&&\ref{1}&$B_1$&$1$&$1$&$1$&$3$& \\
\hline
\cite[Thm.~4.1]{ZHF}&&\ref{1}&$B_2$&$1$&$1$&$1$&$3$& \\
\hline
\cite[Conj.~1(1)]{TZHL}&&\ref{1}&$B_1$&$1$&$1$&$1$&$q+4$&$q\equiv 2\pmod 3$ \\
\hline
\cite[Conj.~1(2)]{TZHL}&&\ref{1}&$B_2$&$1$&$1$&$1$&$q^2-q+1$&$q\equiv 2\pmod 3$ \\
\hline
\cite[Thm.~3.2]{DQWYY}&&\ref{1}&$B_1$&$1$&$1$&$1$&$q^2-q+1$& \\
\hline
\cite[Thm.~4.9]{LQC}&s&\ref{1}&$B_1$&$1$&$1$&$1$&$q^2-q+1$& \\
\hline
\cite[Thm.~3.6]{LQCL}&s&\ref{1}&$B_1$&$1$&$1$&$1$&$q^2-q+1$& \\
\hline

\cite[Thm.~3.15]{WYDM}&&\ref{1}&$B_1$&$1$&$1$&$1$&& \\
\hline

\cite[Thm.~3.4]{GS}&&\ref{1}&$B_2$&$1$&$2$&$2$&$5$& \\
\hline

\cite[Thm.~3.8]{LQCL}&m&\ref{1}&$B_2$&$1$&$2$&$2$&$q^2-2q+2$&$q\equiv 2\pmod 3$ \\
\hline
\cite[Thm.~6]{LH-sev}&m&\ref{1}&$B_2$&$1$&$2$&$2$&$q^2-2q+2$& \\
\hline
\cite[Thm.~3.5]{GS}&m&\ref{1}&$B_2$&$1$&$2$&$2$&$q^2-q+3$& \\
\hline
\cite[Thm.~2.8]{LQLF}&&\ref{1}&$B_2$&$1$&$2$&$2$&& \\
\hline
\cite[Thm.~4.4]{ZHF}&&\ref{2}&$B_1$&$2$&$1$&$2$&$5$& \\
\hline
\cite[Conj.~2]{GS}&&\ref{2}&$B_2$&$2$&$1$&$2$&$5$& \\
\hline
\cite[Thm.~4.3]{ZHF}&&\ref{2}&$B_2$&$2$&$1$&$2$&$5$& \\
\hline
\cite[Thm.~3]{LH-sev}&m&\ref{2}&$B_2$&$2$&$1$&$2$&$q^2-2q+2$& \\
\hline
\cite[Cor.~4.3]{WYDM}&&\ref{2}&$B_2$&$2$&$1$&$2$&$q^2-2q+2$& \\
\hline
\cite[Thm.~4.4]{WYDM}&&\ref{2}&$B_1$&$2$&$1$&$2$&& \\
\hline
\cite[Thm.~2.4]{LQLF}&&\ref{2}&$B_2$&$2$&$1$&$2$&& \\
\hline
\cite[Thm.~4.2]{WYDM}&&\ref{2}&$B_2$&$2$&$1$&$2$&& \\
\hline

\cite[Thm.~3.3]{SHK}&&\ref{1}&$B_2$&$2$&$2$&$2$&$6$& \\
\hline
\cite[Thm.~3.4]{SHK}&m&\ref{1}&$B_2$&$2$&$2$&$2$&$6$& \\
\hline
\end{tabular}
\end{table}

\begin{table}[!htbp]
\setlength\extrarowheight{3pt}
\caption{Previous results and conjectures subsumed by our results, II}
\label{tabpre2}
\begin{tabular}{l|c|c|c|c|c|c}
\hline
Reference&Thm.&$B_z$&$Q$&$R$&$S$&$r$ \\
\hline
\hline
\cite[Thm.~3.5]{XCP}&\ref{1}&$B_1$&$4$&$1$&$2$&$7$ \\
\hline
\cite[Thm.~3.2]{ZKZPZ}&\ref{1}&$B_1$&$4$&$1$&$2$&$7$ \\
\hline
\cite[Thm.~4.2]{XCP}&\ref{2}&$B_2$&$4$&$1$&$2$&$7$ \\
\hline
\cite[Thm.~3.6]{XCP}&\ref{1}&$B_2$&$4$&$1$&$2$&$q^2-q+5$ \\
\hline
\cite[Thm.~3.5]{ZKZPZ}&\ref{1}&$B_2$&$4$&$1$&$2$&$q^2-q+5$ \\
\hline
\cite[Thm.~4.9]{XCP}&\ref{2}&$B_2$&$4$&$1$&$4$&$q^2-q+7$ \\
\hline
\cite[Thm.~3.6]{ZKZPZ}&\ref{2}&$B_2$&$4$&$1$&$4$&$q^2-q+7$ \\
\hline
\cite[Thm.~3.3]{ZZZ}&\ref{1}&$B_2$&$1$&$8$&$2$&$11$ \\
\hline
\cite[Thm.~4.8]{ZZZ}&\ref{1}&$B_1$&$1$&$8$&$2$& \\
\hline
\cite[Thm.~3.6]{ZZZ}&\ref{1}&$B_1$&$4$&$1$&$8$&$13$ \\
\hline
\cite[Thm.~4.2]{ZZZ}&\ref{2}& $B_2$&$4$&$1$&$8$&$13$ \\
\hline
\cite[Table~2(5)]{ZZZ}&\ref{1}&$B_1$&$16$&$1$&$2$&$19$ \\
\hline
\cite[Thm.~3.8]{ZZZ}&\ref{1}&$B_1$&$16$&$4$&$1$&$21$\\
\hline
\cite[Thm.~3.11]{ZZZ}&\ref{1}&$B_1$&$16$&$1$&$8$&$25$ \\
\hline
\cite[Table~2(9)]{ZZZ}&\ref{1}&$B_2$&$1$&$32$&$2$&$35$ \\
\hline
\cite[Thm.~3.13]{ZKZPZ}&\ref{1}&$B_2$&$1$&$32$&$2$&$35$ \\
\hline
\cite[Table~2(10)]{ZZZ}&\ref{1} &$B_1$&$4$&$1$&$32$&$37$ \\
\hline
\cite[Thm.~4.6]{ZZZ}&\ref{2}&$B_1$&$1$&$8$&$32$&$41$ \\
\hline
\cite[Table~2(11)]{ZZZ}&\ref{1}& $B_2$&$1$&$32$&$8$&$41$ \\
\hline
\cite[Table~2(12)]{ZZZ}&\ref{1}&$B_1$&$16$&$1$&$32$&$49$ \\
\hline
\cite[Thm.~3.10]{ZKZPZ}&\ref{1}&$B_1$&$16$&$1$&$32$&$49$ \\
\hline
\cite[Table~2(13)]{ZZZ}&\ref{1} &$B_1$&$64$&$1$&$2$&$67$ \\
\hline
\cite[Table~2(14)]{ZZZ}&\ref{1}&$B_1$&$64$&$4$&$1$&$69$ \\
\hline
\cite[Thm.~3.18]{ZKZPZ}&\ref{1}&$B_1$&$64$&$4$&$1$&$69$ \\
\hline
\cite[Table~2(15)]{ZZZ}&\ref{1}&$B_1$&$64$&$1$&$8$&$73$ \\
\hline
\cite[Table~2(16)]{ZZZ}&\ref{1}&$B_1$&$64$&$16$&$1$&$81$ \\
\hline
\cite[Table~2(17)]{ZZZ}&\ref{1}&$B_1$&$64$&$1$&$32$& $97$ \\
\hline
\end{tabular}
\end{table}

\begin{table}[!htbp]
\setlength\extrarowheight{3pt}
\caption{Previous results and conjectures subsumed by our results, III}
\label{tabpre3}
\begin{tabular}{l|c|c|c|c|c|c|c|l}
\hline
Reference&e&Thm.&$B_z$&$Q$&$R$&$S$&$r$&Conditions \\
\hline
\hline

\multirow{2}{*}{\cite[Ex.~1(4)]{DZ}}&&\ref{2}&$B_2$&&$1$&$2$&$Q+3$&$q\equiv Q\equiv 1\pmod 3$ \\
&&\ref{2}&$B_1$&$1$&$2$&&$S+3$&$S\equiv 2\pmod 3$ \\
\hline
\multirow{2}{*}{\cite[Ex.~1(5)]{DZ}}&&\ref{2}&$B_2$&&$2$&$1$&$Q+3$&$Q\equiv 1\pmod 3$ \\
&&\ref{2}&$B_1$&&$1$&$2$&$Q+3$&$Q\equiv 2\pmod 3$ \\
\hline

\cite[Thm.~2.3]{SV}&s&\ref{1}&$B_2$&$1$&&$R$&$2R+1$&$q\equiv R\equiv 2\pmod 3$ \\
\hline
\cite[Thm.~2.4]{SV}&s&\ref{1}&$B_1$&$1$&&$R$&$2R+1$&$q\equiv R\equiv 2\pmod 3$ \\
\hline
\cite[Thm.~2.5]{SV}&s&\ref{1}&$B_1$&&$Q$&$1$&$2Q+1$&$Q\equiv 1\pmod 3$ \\
\hline
\cite[Thm.~2.6]{SV}&s&\ref{1}&$B_2$&&$Q$&$1$&$2Q+1$&$Q\equiv 1\pmod 3$ \\
\hline
\multirow{2}{*}{\cite[Thm.~1]{LH-new}}&m&\ref{1}&$B_1$&&$Q$&$1$&$Q(q^2-q)+1$&$Q\equiv 1\pmod 3$ \\
&m&\ref{2}&$B_2$&&$1$&$Q$&$Q(q^2-q)+1$&$Q\equiv 2\pmod 3$ \\
\hline
\multirow{2}{*}{\cite[Thm.~2]{LH-new}}&m&\ref{1}&$B_2$&$1$&&$R$&$R(q^2-q)+1$&$R\equiv 2\pmod 3$ \\
&m&\ref{2}&$B_1$&&$1$&$Q$&$Q(q^2-q)+1$& \\
\hline
\cite[Thm.~3.1]{WZZ}&m&\ref{2}&$B_2$&&$1$&$Q$&$Q(q^2-q)+1$&$Q\equiv 1\pmod 3$ \\
\hline
\multirow{3}{*}{\cite[Thm.~2]{KX}}&&\ref{2}&$B_1$&&$1$&&$Q+S+1$&$Q\equiv S\equiv 1\pmod 3$ \\
&&\ref{2}&$B_2$&$1$&&&$R+S+1$&$S\not\equiv R\equiv 1\pmod 3$ \\
&&\ref{1}&$B_2$&$1$&&&$R+S+1$&$R\equiv S\equiv 2\pmod 3$ \\
\hline
\multirow{4}{*}{\cite[Thm.~3]{KX}}&&\ref{2}&$B_1$&&&$1$&$Q+R+1$&$Q\equiv R\equiv 1\pmod 3$ \\
&&\ref{2}&$B_2$&&$1$&&$Q+S+1$&$S\not\equiv Q\equiv 1\pmod 3$ \\
&&\ref{2}&$B_1$&$1$&&&$R+S+1$&$R\equiv S\equiv 1\pmod 3$ \\
&&\ref{1}&$B_1$&&$1$&&$Q+S+1$&$S\not\equiv Q\equiv 1\pmod 3$ \\
\hline
\multirow{3}{*}{\cite[Thm.~4]{KX}}&&\ref{2}&$B_1$&&&$1$&$Q+R+1$&$R\not\equiv Q\equiv 1\pmod 3$ \\
&&\ref{2}&$B_2$&&$1$&&$Q+S+1$&$Q\equiv S\equiv 2\pmod 3$ \\
&&\ref{1}&$B_1$&&&$1$&$Q+R+1$&$Q\equiv R\equiv 1\pmod 3$ \\
\hline
\cite[Thm.~3.1]{ZKZPZ}&&\ref{1}&$B_1$&&$1$&&$Q+S+1$&$S\not\equiv Q\equiv 1\pmod 3$ \\
\hline
\cite[Thm.~3.8]{ZKZPZ}&&\ref{1}&$B_1$&&$1$&&$Q+S+1$&$S\not\equiv Q\equiv 1\pmod 3$ \\
\hline
\cite[Thm.~3.12]{ZKZPZ}&&\ref{1}&$B_1$&$1$&&&$R+S+1$&$R\equiv S\equiv 2\pmod 3$ \\
\hline
\cite[Thm.~3.17]{ZKZPZ}&&\ref{1}&$B_1$&&&$1$&$Q+R+1$&$Q\equiv R\equiv 1\pmod 3$ \\
\hline

\end{tabular}
\end{table}

\begin{table}[!htbp]
\setlength\extrarowheight{3pt}
\caption{Previous results and conjectures subsumed by our results, IV}
\label{tabpre4}
\begin{tabular}{l|c|c|c|c|c|c|l}
\hline
Reference&Thm.&$B_z$&$Q$&$R$&$S$&$r$&Conditions \\
\hline
\hline

\cite[Cor.~3.7]{WYDM}&\ref{1}&$B_1$&&$Q$&&&$S\not\equiv Q\equiv q\equiv 1\pmod 3$ \\
\hline
\cite[Cor.~3.8]{WYDM}&\ref{1}&$B_1$&&$Q$&&&$q\equiv Q\equiv S\equiv 1\pmod 3$ \\
\hline
\cite[Cor.~3.9]{WYDM}&\ref{1}&$B_1$&&$Q$&&&$q\equiv Q\equiv S\equiv 1\pmod 3$ \\
\hline
\cite[Cor.~3.12]{WYDM}&\ref{1}&$B_1$&&$Q$&&&$q\not\equiv Q\equiv S\equiv 1\pmod 3$ \\
\hline
\cite[Cor.~3.13]{WYDM}&\ref{1}&$B_1$&&$Q$&&&$q\not\equiv Q\equiv S\equiv 1\pmod 3$ \\
\hline
\cite[Cor.~3.14]{WYDM}&\ref{1}&$B_2$&&$Q$&&&$Q\not\equiv S\equiv q\equiv 2\pmod 3$ \\
\hline
\cite[Cor.~3.10]{WYDM}&\ref{1}&$B_1$&&&$R$&&$R\not\equiv Q\equiv q\equiv 1\pmod 3$ \\
\hline
\cite[Cor.~3.11]{WYDM}&\ref{1}&$B_2$&&&$R$&&$Q\not\equiv R\equiv q\equiv 2\pmod 3$ \\
\hline
\multirow{4}{*}{\cite[Thm.~3.1]{GWSHZL}}&\ref{1}&$B_1$&&$Q$&&&$Q\equiv S\equiv 1\pmod 3$ \\
&\ref{1}&$B_2$&&$Q$&&&$Q\equiv S\equiv 2\pmod 3$ \\
&\ref{2}&$B_1$&&&$Q$&&$R\not\equiv Q\equiv 1\pmod 3$ \\
&\ref{2}&$B_2$&&&$Q$&&$Q\not\equiv R\equiv 1\pmod 3$ \\
\hline
\multirow{4}{*}{\cite[Thm.~1.1]{Z-x}}&\ref{1}&$B_1$&&$Q$&&&$Q\equiv S\equiv 1\pmod 3$ \\
&\ref{1}&$B_2$&&$Q$&&&$Q\equiv S\equiv 2\pmod 3$ \\
&\ref{2}&$B_1$&&&$Q$&&$R\not\equiv Q\equiv 1\pmod 3$ \\
&\ref{2}&$B_2$&&&$Q$&&$Q\not\equiv R\equiv 1\pmod 3$ \\
\hline
\end{tabular}
\end{table}

\FloatBarrier

In our final result, we show that in the special case $r=Q+R+S$ the polynomials $f(X)$ in Theorems~\ref{1} and \ref{2} can be expressed as compositions of functions having very simple forms. The result relies on the following notion.

\begin{defn}
If $U$ and $V$ are $\F_q$-vector spaces, then a function $f\colon U\to U$ is \emph{$\F_q$-linearly equivalent} to a function $g\colon V\to V$ if $f=\rho\circ g\circ\eta$ for some $\F_q$-vector space isomorphisms $\rho\colon V \to U$ and $\eta\colon U\to V$.
\end{defn}

\begin{rmk}
It is easy to see that $\F_q$-linear equivalence is an equivalence relation on the union of the sets of functions $\F_{q^2}\to\F_{q^2}$ and $\F_q\times\F_q\to\F_q\times\F_q$, and that $\F_q$-linear equivalence preserves the property of a function being bijective.
\end{rmk}

\begin{thm}\label{3}
In the notation of Theorem~\emph{\ref{1}} or \emph{\ref{2}}, we have:
\begin{itemize}
\item if $r=Q+R+S$ and $q\equiv 1\pmod 3$ then $f(X)$ in Theorem~\emph{\ref{1}} is\/ $\F_q$-linearly equivalent to $X^{Q+R+S}$ on\/ $\F_{q^2}$; 
\item if $r=Q+R+S$ and $q\equiv 2\pmod 3$ then $f(X)$ in Theorem~\emph{\ref{1}} is\/ $\F_q$-linearly equivalent to $(X^{Q+R+S},Y^{Q+R+S})$ on\/ $\F_q\times\F_q$; 
\item if $r=Q+R+S$ and $q\equiv 1\pmod 3$ then $f(X)$ in Theorem~\emph{\ref{2}} is\/ $\F_q$-linearly equivalent to $X^{Q+qR+S}$ on\/ $\F_{q^2}$;
\item if $r=Q+R+S$ and $q\equiv 2\pmod 3$ then $f(X)$ in Theorem~\emph{\ref{2}} is\/ $\F_q$-linearly equivalent to $(X^{Q+S}Y^R,X^RY^{Q+S})$ on\/ $\F_q\times\F_q$. 
\end{itemize}
\end{thm}

\begin{rmk}
In light of the well-known fact that $X^n$ permutes $\F_{q^i}$ if and only if $n>0$ and $\gcd(n,q^i-1)=1$, Theorem~\ref{3} gives an alternative proof of the special case $r=Q+R+S$ of Theorems~\ref{1} and \ref{2}.
Although our proof of Theorem~\ref{3} takes only roughly one-third of one page, already the special case $\charp(\F_q)=2$ of Theorem~\ref{3} subsumes many previous results, including for instance all results from the very recent 23-page and 31-page papers \cite{KX} and \cite{ZKZPZ} (cf.\ Table~\ref{tabpre3}).
\end{rmk}

This paper is organized as follows. In the next section we provide background material. Then we prove Theorems~\ref{1}, \ref{2}, and \ref{3} in Sections~\ref{sec1}, \ref{sec2}, and \ref{sec3}, respectively.


\section{Background results}

In this section we present the known results which are used in our proof of Theorems~\ref{1} and ~\ref{2}. They rely on the following notation.

\begin{notation}
If $q$ is a prime power, then we write $\mu_{q+1}$ for the set of all $(q+1)$-th roots of unity in $\F_{q^2}$, and we define $\bP^1(\F_q)\colonequals \F_q\cup\{\infty\}$.
\end{notation}

We begin with the following special case of \cite[Lemma~2.1]{Zlem}.

\begin{lemma}\label{old}
Write $f(X)\colonequals X^r B(X^{q-1})$ where $q$ is a prime power, $r$ is a positive integer, and $B(X)\in\F_{q^2}[X]$. Then $f(X)$ permutes\/ $\F_{q^2}$ if and only if $\gcd(r,q-1)=1$ and $X^r B(X)^{q-1}$ permutes $\mu_{q+1}$.
\end{lemma}

The next two results are reformulations of \cite[Lemmas~2.1 and 3.1]{Z-Redei}.

\begin{lemma}\label{deg1mu}
If $\alpha,\beta\in\F_{q^2}$ satisfy $\alpha^{q+1}\ne\beta^{q+1}$, then $(\beta^qX+\alpha^q)/(\alpha X+\beta)$ 
permutes $\mu_{q+1}$.
\end{lemma}

\begin{lemma}\label{mu}
If $\alpha\in\F_{q^2}\setminus\F_q$ and $\beta\in\mu_{q+1}$, then $(\alpha X+\beta\alpha^q)/(X+\beta)$ maps $\mu_{q+1}$ bijectively onto\/ $\bP^1(\F_q)$.
\end{lemma}

\begin{prop}\label{cubic}
Assume $q$ is a prime power with $\gcd(3,q)=1$. Pick an order-$3$ element $\omega\in\F_{q^2}^*$. Write $\rho(X) \colonequals (X-\omega)/(-\omega X+1)$ and  $\eta(X) \colonequals (X+\omega)/(\omega X+1)$. Then the following hold:
\begin{enumerate}
\item if $q\equiv 1\pmod 3$ then $\rho(X)$ and $\eta(X)$ permute $\mu_{q+1}$; \item if $q\equiv 2\pmod 3$ then $\rho(X)$ and $\eta(X)$ interchange $\mu_{q+1}$ with\/ $\bP^1(\F_q)$.  
\end{enumerate}
\end{prop}

\begin{proof}
If $q\equiv 1\pmod 3$ then $(-\omega)^{q+1}=\omega^{q+1}\ne 1$, which by Lemma~\ref{deg1mu} implies that $\rho(X)$ and $\eta(X)$ permute $\mu_{q+1}$. 
If $q\equiv 2\pmod 3$ then $\omega\in\F_{q^2}\setminus\F_q$ and $\omega\in\mu_{q+1}$, which by Lemma~\ref{mu} implies that $1/\rho(X)$ and $1/\eta(X)$ induce bijections $\mu_{q+1}\to\bP^1(\F_q)$, so that also $\rho(X)$ and $\eta(X)$ induce bijections $\mu_{q+1}\to\bP^1(\F_q)$. 
Since $\rho(X)$ is the inverse of $\eta(X)$ under composition, it follows that $\rho(X)$ and $\eta(X)$ also induce bijections $\bP^1(\F_q)\to\mu_{q+1}$. This concludes the proof.
\end{proof}


\section{Proof of Theorem~\ref{1}}
\label{sec1}

We now prove Theorem~\ref{1}.
Define $\rho(X) \colonequals (X-\omega)/(-\omega X+1)$ and $\eta(X) \colonequals (X+\omega)/(\omega X+1)$.
Then
\[
X^{Q+R+S} \circ \eta(X) =
\frac{(X+\omega)^{Q+R+S}}{(\omega X+1)^{Q+R+S}} =
\frac{N(X)}{D(X)}
\]
where
\begin{align*}
N(X)&\colonequals (X+\omega)^{Q+R+S} \\
&= (X+\omega)^Q (X+\omega)^R (X+\omega)^S \\
&= (X^Q+\omega^Q)(X^R+\omega^R)(X^S+\omega^S)
\end{align*}
and
\begin{align*}
D(X)&\colonequals (\omega X+1)^{Q+R+S} \\
&= (\omega^Q X^Q+1)(\omega^R X^R+1)(\omega^S X^S+1).
\end{align*}
Since $X+\omega$ and $\omega X+1$ are coprime, also $N(X)$ and $D(X)$ are coprime. Define
\[
g(X)\colonequals \rho(X)\circ X^{Q+R+S}\circ\eta(X),
\]
so that $g(X)=U(X)/V(X)$ where $U(X)\colonequals N(X)-\omega D(X)$ and $V(X)\colonequals -\omega N(X)+D(X)$. Here
\[
\gcd\bigl(U(X),V(X)\bigr)=\gcd\bigl(N(X),D(X)\bigr)=1.
\]
One sees by inspection that $U(X)=C_2(X)$ and $V(X)=C_1(X)$. Permuting the values of $Q,R,S$ does not change $U(X)$ or $V(X)$, and has the effect of multiplying $\beta$ by some $\epsilon\in\{1,-1\}$, so it has the effect of multiplying each $B_z(X)$ by $\epsilon$. 
Thus if $a,b,c\in\{1,-1\}$ satisfy $Q\equiv a\pmod 3$, $R\equiv b\pmod 3$, and $S\equiv c\pmod 3$, then there is no loss in assuming $a\ge b\ge c$, and in this case one sees by inspection that $B_z(X)$ is given by the expression in Table~\ref{tab1}. In particular, the coefficients of $B_z(X)$ are in $\F_p$.

By Lemma~\ref{old}, $f(X)\colonequals X^r B_z(X^{q-1})$ permutes $\F_{q^2}$ if and only if $\gcd(r,q-1)=1$ and $X^r B_z(X)^{q-1}$ permutes $\mu_{q+1}$. Note that $B_{3-z}(X)$ equals $X^{Q+R+S} B_z(1/X)$, which induces the same function on $\mu_{q+1}$ as does $X^{Q+R+S} B_z(X^q)$, which equals $X^{Q+R+S} B_z(X)^q$ since $B_z(X)\in\F_q[X]$. 
It follows that any root of $B_z(X)$ in $\mu_{q+1}$ would also be a root of $B_{3-z}(X)$, which is impossible since $B_1(X)$ and $B_2(X)$ are coprime. Thus $B_z(X)$ has no roots in $\mu_{q+1}$. Moreover, $X^r B_z(X)^{q-1}$ induces the same function on $\mu_{q+1}$ as does $X^{Q+R+S} B_z(X)^q/B_z(X)$, which induces the same function on $\mu_{q+1}$ as does $B_{3-z}(X)/B_z(X)$. 
Since $B_{3-z}(X)/B_z(X)$ equals either $g(X)$ or $1/g(X)$, we see that $X^r B_z(X)^{q-1}$ permutes $\mu_{q+1}$ if and only if $g(X)$ permutes $\mu_{q+1}$. If $q\equiv 1\pmod 3$ then $\rho(X)$ and $\eta(X)$ permute $\mu_{q+1}$ by Proposition~\ref{cubic}, so that $g(X)$ permutes $\mu_{q+1}$ if and only if $X^{Q+R+S}$ permutes $\mu_{q+1}$, or equivalently $\gcd(Q-R+S,q+1)=1$. 
If $q\equiv 2\pmod 3$ then $\rho(X)$ and $\eta(X)$ interchange $\mu_{q+1}$ and $\bP^1(\F_q)$ by Proposition~\ref{cubic}, so that $g(X)$ permutes $\mu_{q+1}$ if and only if $X^{Q+R+S}$ permutes $\bP^1(\F_q)$, or equivalently $\gcd(Q+R+S,q-1)=1$.  This concludes the proof.


\section{Proof of Theorem~\ref{2}}
\label{sec2}

We now prove Theorem~\ref{2}.
Define $\rho(X) \colonequals (X-\omega)/(-\omega X+1)$ and $\eta(X) \colonequals (X+\omega)/(\omega X+1)$.
Then
\[
X^{Q-R+S}\circ\eta(X) = \frac{(X+\omega)^{Q-R+S}}{(\omega X+1)^{Q-R+S}}=
\frac{(X+\omega)^{Q+S}(\omega X+1)^R}{(\omega X+1)^{Q+S}(X+\omega)^R} =
\frac{N(X)}{D(X)}
\]
where
\begin{align*}
N(X)&\colonequals 
(X+\omega)^{Q+S}(\omega X+1)^R \\
&=
(X+\omega)^Q(\omega X+1)^R(X+\omega)^S \\
&= (X^Q+\omega^Q)(\omega^R X^R+1)(X^S+\omega^S)
\end{align*}
and
\begin{align*}
D(X)&\colonequals  
(\omega X+1)^{Q+S} (X+\omega)^R \\
&=
(\omega X+1)^Q (X+\omega)^R (\omega X+1)^S \\
&=
(\omega^Q X^Q+1)(X^R+\omega^R)(\omega^S X^S+1).
\end{align*}
Since $X+\omega$ and $\omega X+1$ are coprime, and their product is $\omega (X^2-X+1)$, we have $\gcd\bigl(N(X),D(X)\bigr)=(X^2-X+1)^{\min(Q+S,R)}$. Define
\[
g(X)\colonequals \rho(X)\circ X^{Q-R+S}\circ\eta(X),
\]
so that $g(X)=U(X)/V(X)$ where $U(X)\colonequals N(X)-\omega D(X)$ and $V(X)\colonequals -\omega N(X)+D(X)$. Here 
\[
\gcd\bigl(U(X),V(X)\bigr)=\gcd\bigl(N(X),D(X)\bigr)=(X^2-X+1)^{\min(Q+S,R)}.
\]
By inspection, we see that $U(X)=C_2(X)$ and $V(X)=C_1(X)$. Interchanging $Q$ and $S$ does not change $U(X)$ or $V(X)$, and has the effect of multiplying $\beta$ by some $\epsilon\in\{1,-1\}$, so it has the effect of multiplying each $B_z(X)$ by $\epsilon$. 
By performing such an interchange if necessary, we may assume that if $Q\equiv -1\pmod 3$ then $S\equiv -1\pmod 3$. By inspection, $B_z(X)$ is given by the expression in Table~\ref{tab2}. In particular, the coefficients of $B_z(X)$ are in $\F_p$.

By Lemma~\ref{old}, $f(X)\colonequals X^r B_z(X^{q-1})$ permutes $\F_{q^2}$ if and only if $\gcd(r,q-1)=1$ and $X^r B_z(X)^{q-1}$ permutes $\mu_{q+1}$. Plainly if $X^r B_z(X)^{q-1}$ permutes $\mu_{q+1}$ then $B_z(X)$ has no roots in $\mu_{q+1}$. 
Note that $B_{3-z}(X)$ equals $X^{Q+R+S} B_z(1/X)$, and induces the same function on $\mu_{q+1}$ as does $X^{Q+R+S} B_z(X^q)$, which equals $X^{Q+R+S} B_z(X)^q$ since $B_z(X)\in\F_q[X]$. It follows that any root of $B_z(X)$ in $\mu_{q+1}$ is also a root of $B_{3-z}(X)$, so that $B_z(X)$ has no roots in $\mu_{q+1}$ if and only if $\gcd\bigl(B_1(X),B_2(X)\bigr)$ has no roots in $\mu_{q+1}$. 
We showed above that this gcd is a power of $X^2-X+1$, so the gcd has no roots in $\mu_{q+1}$ if and only if $q\equiv 1\pmod 3$. Henceforth assume $q\equiv 1\pmod 3$. Then $X^r B_z(X)^{q-1}$ induces the same function on $\mu_{q+1}$ as does $X^{Q+R+S} B_z(X)^q/B_z(X)$, which induces the same function on $\mu_{q+1}$ as does $B_{3-z}(X)/B_z(X)$. 
Since $B_{3-z}(X)/B_z(X)$ equals either $g(X)$ or $1/g(X)$, we see that $X^r B_z(X)^{q-1}$ permutes $\mu_{q+1}$ if and only if $g(X)$ permutes $\mu_{q+1}$. 
Since $q\equiv 1\pmod 3$, Proposition~\ref{cubic} implies that $\rho(X)$ and $\eta(X)$ permute $\mu_{q+1}$, so that $g(X)=\rho(X)\circ X^{Q-R+S}\circ\eta(X)$ permutes $\mu_{q+1}$ if and only if $X^{Q-R+S}$ permutes $\mu_{q+1}$, or equivalently $\gcd(Q-R+S,q+1)=1$. This concludes the proof.


\section{Proof of Theorem~\ref{3}}\label{sec3}

We now prove Theorem~\ref{3}. Since the permutation property of the polynomials in this result was already resolved in Theorems~\ref{1} and \ref{2}, we decided not to write out the details of some routine verifications in the proof of Theorem~\ref{3}, which any interested reader could easily check.

First suppose $q\equiv 1\pmod 3$. Let $\eta$ be the $\F_q$-vector space automorphism of $\F_{q^2}$ induced by $\omega X^q+X$. Let $\rho$ be the $\F_q$-vector space automorphism of $\F_{q^2}$ which is induced by $\beta^{-1}(-\omega X^q+X)$ if $z=1$, and by $\beta^{-1}(X^q-\omega X)$ if $z=2$. 
Let $g\colon\F_{q^2}\to\F_{q^2}$ be the function induced by $X^{Q+R+S}$ (resp., $X^{Q+qR+S}$) for Theorem~\ref{1} (resp., Theorem~\ref{2}). It is routine to verify that the function on $\F_{q^2}$ induced by $f(X)$ equals $\rho\circ g\circ\eta$.

Next suppose $q\equiv 2\pmod 3$. The function $\eta\colon x\mapsto \bigl(\omega^q x+(\omega^q x)^q,\omega x+(\omega x)^q\bigr)$ induces an $\F_q$-vector space isomorphism $\F_{q^2}\to\F_q\times\F_q$.
Let $\rho\colon\F_q\times\F_q\to\F_{q^2}$ be the $\F_q$-vector space isomorphism  which is defined by $(x,y)\mapsto \beta^{-1}\omega^{-Q-R-S}(-\omega x+y)$ if $z=1$, and by $(x,y)\mapsto \beta^{-1}\omega^{-Q-R-S}(x-\omega y)$ if $z=2$.
Let $g$ be the function $\F_q\times\F_q\to\F_q\times\F_q$ induced by $(X^{Q+R+S},Y^{Q+R+S})$ (resp., $(X^{Q+S}Y^R,X^RY^{Q+S})$) for Theorem~\ref{1} (resp., Theorem~\ref{2}). It is routine to verify that the function on $\F_{q^2}$ induced by $f(X)$ equals $\rho\circ g\circ\eta$.

\end{document}